\documentclass[12pt]{article}

\usepackage{amsfonts,graphicx,color}

\title{Compact hyperbolic tetrahedra with non-obtuse dihedral angles}

\author{Roland K. W. Roeder \footnote{rroeder@fields.utoronto.ca}}
                                                                                    
\input{psfig.sty}
                                                                                    
\newtheorem{thm}{Theorem}
\newtheorem{prop}[thm]{Proposition}
\newtheorem{cor}[thm]{Corollary}
\newtheorem{lem}[thm]{Lemma}
                                                                                    
\newcommand{\Endproof}{$\Box$ \vspace{.15in}}

\begin{document}

\maketitle

\begin{abstract}

Given a combinatorial description $C$ of a polyhedron having $E$ edges, the
space of dihedral angles of all compact hyperbolic polyhedra that realize $C$
is generally not a convex subset of $\mathbb{R}^E$ \cite{DIAZ}.  If $C$ has
five or more faces, Andreev's Theorem states that the corresponding space of
dihedral angles $A_C$ obtained by restricting to {\em non-obtuse} angles is a
convex polytope.  In this paper we explain why Andreev did not consider
tetrahedra, the only polyhedra having fewer than five faces, by demonstrating
that the space of dihedral angles of compact hyperbolic tetrahedra, after
restricting to non-obtuse angles, is non-convex.  Our proof provides a simple
example of the ``method of continuity'', the technique used in classification
theorems on polyhedra by Alexandrow \cite{ALEX}, Andreev \cite{AND}, and
Rivin-Hodgson \cite{RH}. \\


\end{abstract}
Given a combinatorial description $C$ of a polyhedron having $E$ edges, 
the space of dihedral angles of all compact hyperbolic polyhedra that 
realize $C$ is generally not a convex subset of $\mathbb{R}^E$.  This is 
proved in a nice paper by D\'iaz \cite{DIAZ}. However, Andreev's Theorem 
\cite{AND,H,ROE2,ROE} shows that by restricting to compact hyperbolic 
polyhedra with non-obtuse dihedral angles, the space of dihedral angles is 
a convex polytope, which we label $A_C \subset \mathbb{R}^E$.  It is 
interesting to note that the statement of Andreev's Theorem requires that 
$C$ have five or more faces, ruling out the tetrahedron which is the only 
polyhedron having fewer than five faces.

In this paper, we explain why hyperbolic tetrahedra are a special case that is
not covered by Andreev's Theorem.  We provide an explicit description of the
space of dihedral angles, $A_\Delta$, corresponding to compact hyperbolic
tetrahedra with non-obtuse dihedral angles, finding that $A_\Delta$ is a
non-convex, path-connected subset of $\mathbb{R}^6$.  

A description of the space of Gram matrices (and hence indirectly of the space
of dihedral angles) corresponding to compact hyperbolic tetrahedra having
arbitrary dihedral angles is available in Milnor's collected works \cite{MIL}.
Our description of the space of dihedral angles $A_\Delta$ can be derived from
the result in \cite{MIL}, using the assumption that the dihedral angles are
non-obtuse.  However, we use the ``method of continuity,'' providing the reader
with a simple example of a method that plays an important role in the
classification theorems on polyhedra by  Alexandrow \cite{ALEX}, Andreev
\cite{AND}, and Rivin-Hodgson \cite{RH}.

\vspace{.1in}

Let $E^{3,1}$ be $\mathbb{R}^4$ with the indefinite metric $\Vert {\bf x} 
\Vert^2 = -x_0^2+x_1^2+x_2^2+x_3^2$. In this paper, we work in the 
hyperbolic space $\mathbb{H}^3$ given by the component of the subset of 
$E^{3,1}$ given by

$$\Vert {\bf x} \Vert^2 = -x_0^2+x_1^2+x_2^2+x_3^2 = -1$$

\noindent
having $x_0 > 0$, with the Riemannian metric induced by the indefinite
metric

$$-dx_0^2+dx_1^2+dx_2^2+dx_3^2.$$

There is a natural compactification of the hyperbolic space obtained by 
adding the set of rays asymptotic to the hyperboloid.  We refer to these 
points as the {\it points at infinity}. There is no natural extension of 
the Riemannian structure of $\mathbb{H}^3$ to these points at infinity, 
however, there is a natural way to extend the conformal structure on 
$\mathbb{H}^3$ to these points at infinity.

\vspace{.05in}

One can check that the hyper-plane orthogonal to a vector ${\bf v} \in 
E^{3,1}$ intersects $\mathbb{H}^3$ if and only if $\langle{\bf v},{\bf 
v}\rangle> 0$.  Let ${\bf v} \in E^{3,1}$ be a vector with $\langle{\bf 
v},{\bf v}\rangle > 0$, and define

$$P_{\bf v} = \{{\bf w} \in \mathbb{H}^3 | \langle{\bf w},{\bf v}\rangle
= 0\}$$

\noindent to be the hyperbolic plane orthogonal to ${\bf v}$;  and the
corresponding closed half space:

$$H_{\bf v}^+ = \{{\bf w} \in \mathbb{H}^3 | \langle{\bf w},{\bf v}\rangle 
\geq 0 \}.$$

\noindent 
Notice that given two planes $P_{\bf v}$ and $P_{\bf w}$ in 
$\mathbb{H}^3$ with $\langle{\bf v},{\bf v}\rangle = 1$ and $\langle{\bf 
w},{\bf w}\rangle = 1$, they:

\begin{itemize}
\item intersect in a line if and only if $\langle{\bf v},{\bf w}\rangle^2
< 1$, in which case their dihedral angle is $\arccos(-\langle{\bf v},{\bf
w}\rangle)$.

\item intersect in a single point at infinity if and only if $\langle{\bf
v},{\bf w}\rangle^2 = 1$, in this case their dihedral angle is $0$.

\end{itemize}

A {\it hyperbolic polyhedron} is an intersection

$$P = \bigcap_{i=0}^n H_{\bf v_i}^+ $$

\noindent having non-empty interior.  There are many papers on hyperbolic
polyhedra, including
\cite{AND,BAO,DIAZ,H,MIL,RH,ROE2,ROE,SCH2,SCH1,VIN,VINVOL}, and particularly on
the groups of reflections generated by them \cite{AVS,VIN,VINREFL,VS}.  A {\it
hyperbolic tetrahedron} is therefore a hyperbolic polyhedron having the
combinatorial type of a tetrahedron.  There are also many papers on hyperbolic
tetrahedra including \cite{CHOKIM,DL,FEL,MARTIN,YANA,SEIDEL}, many of these
studying  volume and symmetries.

If we normalize the vectors $v_i$ that are orthogonal to the faces of a
polyhedron $P$, the {\em Gram Matrix} of $P$ is the matrix with terms
$M_{ij}=v_i \cdot v_j$.  By construction, a Gram matrix is symmetric and
unidiagonal (i.e. has $1$s on the diagonal).  The following Theorem appears in
\cite{MIL}:

\begin{thm}\label{GRAM} 
A symmetric unidiagonal matrix $M$ is the Gram matrix of a compact
hyperbolic tetrahedron if and only if $\det(M) <0$ and each principal
minor is positive definite. 
\end{thm}

\vspace{.1in}

Although the hyperboloid model of hyperbolic space is very natural, it is 
not easy to visualize, since the ambient space is four-dimensional.  We 
will often use the Poincar\'e ball model of hyperbolic space, given by the 
unit ball in $\mathbb{R}^3$ with the metric

$$4\frac{dx_1^2+dx_2^2+dx_3^2}{(1 -\Vert {\bf x}\Vert^2)^2}$$

\noindent
and the upper half-space model of hyperbolic space, given by the subset of 
$\mathbb{R}^3$ with $x_3 > 0$ equipped with the metric 
$$\frac{dx_1^2+dx_2^2+dx_3^2}{x_3^2}.$$

\noindent
Both of these models are isometric to $\mathbb{H}^3$.  The points at 
infinity in the Poincar\'e Ball model correspond to points on the unit 
sphere, and the points at infinity in the upper half-space model 
correspond to the points in the plane $x_3=0$. More background is 
available on hyperbolic geometry in \cite{BEN}.

Hyperbolic planes in these models correspond to portions of Euclidean 
spheres and Euclidean planes that intersect the boundary perpendicularly.  
Furthermore, these models are conformally correct, that is, the hyperbolic 
angle between a pair of such intersecting hyperbolic planes is exactly the 
Euclidean angle between the corresponding spheres or planes.

See below for an image of a compact hyperbolic tetrahedron depicted in the
Poincar\'e ball model depicted using Geomview \cite{GEO}.  The sphere at
infinity is shown for reference.

\begin{center}
\includegraphics[scale = .5]{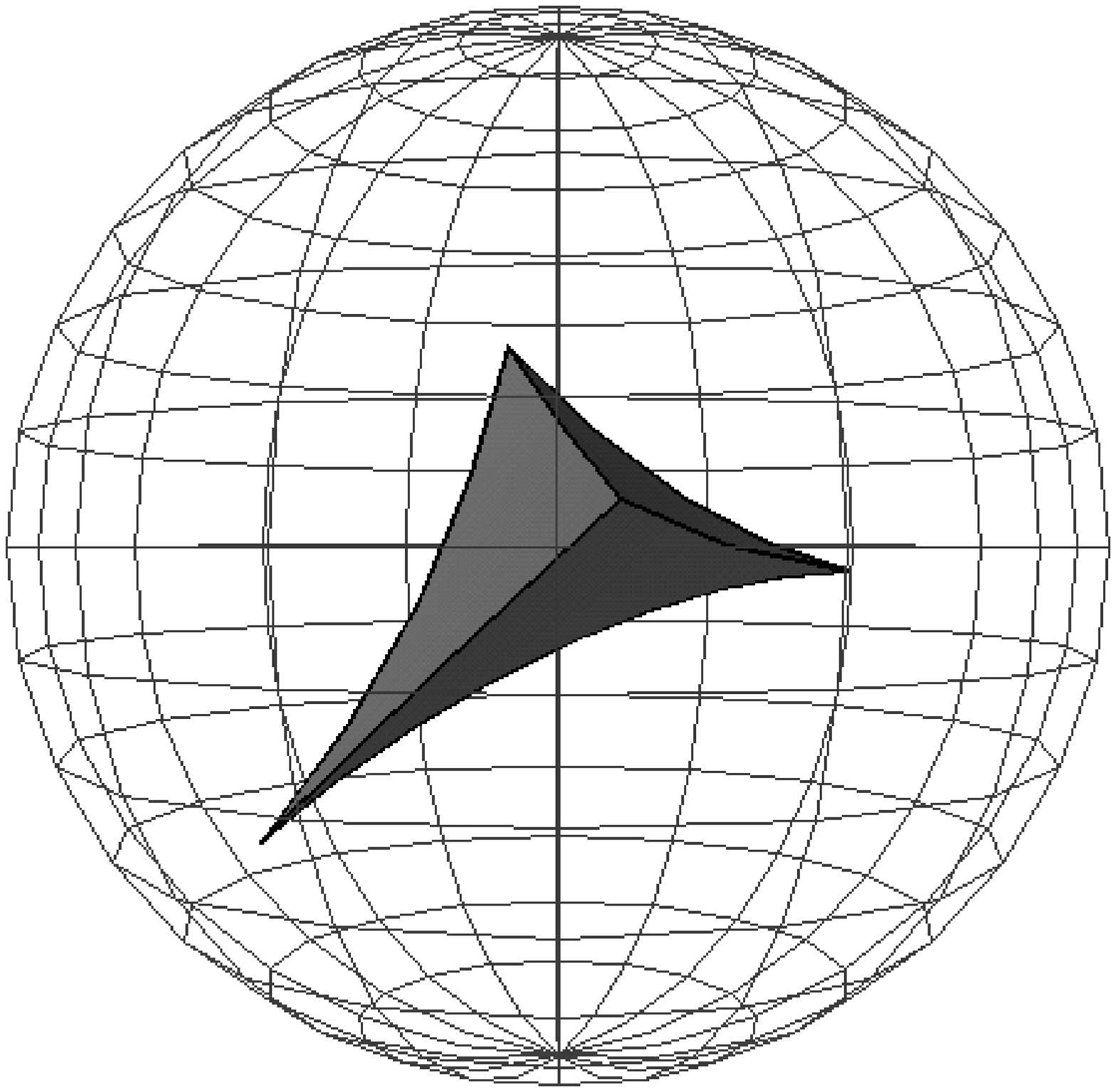}
\end{center}

\noindent

The following two lemmas will be necessary when discussing compact hyperbolic
polyhedra having non-obtuse dihedral angles.  They are well known results and
appear in many of the works on hyperbolic polyhedra mentioned above, including
\cite{AND}.

\begin{lem} \label{INTERSECT}
Suppose that three planes $P_{\bf v_1},P_{\bf v_2},P_{\bf v_3}$ intersect 
pairwise in $\mathbb{H}^3$ with non-obtuse dihedral angles $\alpha, 
\beta$, and $\gamma$.  Then, $P_{\bf v_1},P_{\bf v_2},P_{\bf v_3}$ 
intersect at a vertex in $\overline{\mathbb{H}^3}$ if and only if 
$\alpha+\beta+\gamma \geq \pi.$ The planes intersect in $\mathbb{H}^3$ if 
and only if the inequality is strict.
\end{lem}

\noindent {\bf Proof:}
The planes intersect in a point of $\overline{\mathbb{H}^3}$ if and only 
if the subspace spanned by ${\bf v_1},{\bf v_2},{\bf v_3}$ is positive 
semi-definite, so that the orthogonal is a negative semi-definite line of 
$E^{3,1}$. If the inner product on this line is negative, the line defines 
a point of intersection with the hyperboloid model.  Otherwise, the inner 
product on the line is zero, this line corresponds to a point in 
$\partial\mathbb{H}^3$, since the line will then lie in the cone to which 
the hyperboloid is asymptotic.

The symmetric
matrix defining the inner product on the span of ${\bf v_1},{\bf v_2}$, and ${\bf v_3}$ is
$$\left[
\begin{array}{ccc}
1 & \langle{\bf v_1},{\bf v_2}\rangle &  \langle{\bf v_1},{\bf v_3}\rangle \\
\langle{\bf v_1},{\bf v_2}\rangle & 1 &  \langle{\bf v_2},{\bf v_3}\rangle \\
 \langle{\bf v_1},{\bf v_3}\rangle &  \langle{\bf v_2},{\bf v_3}\rangle & 1\\
\end{array}
\right]
=
\left[
\begin{array}{ccc}
1 & -\cos\alpha & -\cos\beta \\
-\cos\alpha & 1 & -\cos\gamma \\
-\cos\beta & -\cos\gamma & 1 \\
\end{array}
\right]
$$
where $\alpha,\beta,$ and $\gamma$ are the dihedral angles between the pairs of
faces $(P_{\bf v_1},P_{\bf v_2})$, $(P_{\bf v_1},P_{\bf v_3}),$ and $(P_{\bf
v_2},P_{\bf v_3})$, respectively.

Since the principal minor is positive definite for $0 < \alpha \leq 
\pi/2$, it is enough to find out when the determinant $$1 
-2\cos\alpha\cos\beta\cos\gamma -\cos^2\alpha-\cos^2\beta-\cos^2\gamma$$ 
is non-negative.

A bit of trigonometric trickery (we used complex exponentials) shows that 
the expression above is equal to 

{\small
\begin{eqnarray}\label{COSEQN} 
-4\cos\left(\frac{\alpha+\beta+\gamma}{2}\right)\cos\left(\frac{\alpha-\beta+\gamma}{2}\right)\cos\left(\frac{\alpha+\beta-\gamma}{2}\right)\cos\left(\frac{-\alpha+\beta+\gamma}{2}\right) 
\end{eqnarray}
}

Let $\delta = \alpha+\beta+\gamma$.  When $\delta < \pi$, (\ref{COSEQN}) 
is strictly negative; when $\delta = \pi$, (\ref{COSEQN}) is clearly zero; 
and when $\delta > \pi$, (\ref{COSEQN}) is strictly positive.  Hence the 
inner product on the space spanned by ${\bf v_1},{\bf v_2},{\bf v_3 }$ is 
positive semidefinite if and only if $\delta \geq \pi$.  It is positive 
definite if and only if $\delta > \pi$.

Therefore, the three planes $P_{\bf v_1},P_{\bf v_2},P_{\bf v_3} \subset 
\mathbb{H}^3$ intersect at a point in $\overline{\mathbb{H}^3}$ if and 
only if they intersect pairwise in $\mathbb{H}^3$ and the sum of the 
dihedral angles $\delta \geq \pi$.  It is also clear that they intersect 
at a finite point if and only if the inequality is strict. \Endproof

\begin{lem} \label{TRIVALENT}
Given a trivalent vertex of a hyperbolic polyhedron, we can compute the 
angles of the faces in terms of the dihedral angles.  If the dihedral 
angles are non-obtuse, these angles are also $\leq \pi/2$.
\end{lem}

\noindent
{\bf Proof:}
Let $v$ be a finite trivalent vertex of $P$.  After an appropriate 
isometry, we can assume that $v$ is the origin in the Poincar\'e ball 
model, so that the faces at $v$ are subsets of Euclidean planes through 
the origin.  A small sphere centered at the origin will intersect $P$ in a 
spherical triangle $Q$ whose angles are the dihedral angles between faces.  
Call these angles $\alpha_1,\alpha_2,\alpha_3$.

The edge lengths of $Q$ are precisely the angles in the faces at the 
origin. Supposing that $Q$ has edge lengths $(\beta_1,\beta_2,\beta_3)$ 
with the edge $\beta_i$ opposite of angle $\alpha_i$ for each $i = 1,2,3$, 
The law of cosines in spherical geometry states that:

\begin{eqnarray} \label{TLC}
\cos(\beta_i) = \frac{\cos(\alpha_i)
+\cos(\alpha_j)\cos(\alpha_k)} {\sin(\alpha_j)\sin(\alpha_k)}.
\end{eqnarray}

\noindent Hence, the face angles are calculable from the dihedral angles. 
They are non-obtuse, since the right-hand side of the equation is positive 
for $\alpha_i,\alpha_j,\alpha_k$ non-obtuse. \Endproof

\vspace{.1in}
We can now state our classification of compact hyperbolic tetrahedra:

\begin{thm} \label{TETRA1}
Let $\alpha_1,\cdots,\alpha_6$ be a set of proposed non-obtuse dihedral 
angles and let 
$\beta_1(\alpha_1,\cdots,\alpha_6),\cdots,\beta_{12}(\alpha_1,\cdots,\alpha_6)$ 
be the face angles given by equation \ref{TLC}, corresponding to these 
proposed dihedral angles.

\vspace{.1in}
\noindent
There is a compact hyperbolic tetrahedron with dihedral angles
$\alpha_1,\cdots,\alpha_6$ if and only if:
                                                                                
\begin{enumerate}
                                                                                
\item  For each edge $e_i$, $0 < \alpha_i \leq \pi/2$.
                                                                                
\item  Whenever 3 distinct edges $e_i, e_j, e_k$ meet at a vertex,
$\alpha_i + \alpha_j+\alpha_k > \pi$.
                                                                                
\item  For each face the sum of the face angles satisfies
$\beta_i+\beta_j+\beta_k < \pi$.
                                                                                
\end{enumerate}
Furthermore this tetrahedron is unique.
\end{thm}
                                                                                
Recall from Lemma \ref{TRIVALENT} that the face angles $\beta_i$ are 
calculable from the dihedral angles $\alpha_i$ and are themselves 
non-obtuse so that condition (3) is a highly non-linear condition on the 
dihedral angles. We will denote the subset of $\mathbb{R}^6$ of dihedral 
angles satisfying conditions (1-3) by $A_\Delta$.

We present a proof of Theorem \ref{TETRA1} using the ``method of continuity'',
the classical method used by Alexandrow \cite{ALEX}, Andreev \cite{AND}, later
by Rivin and Hodgson \cite{RH}, and in this author's more recent proof of
Andreev's Theorem \cite{ROE2,ROE}.  The idea of this method is to establish a
bijection between two manifolds of the same dimension: one, $X$, consisting of
the geometric objects that you want to construct, and the other, $Y$, a subset
of $\mathbb{R}^n$ consisting of various angles, lengths, etc.  The space $X$
should be viewed as unknown and the space $Y$ as known.

You then consider the mapping $f: X \rightarrow Y$ which takes your
geometric object, in $X$, and reads off its appropriate measurements, in
$Y$.  Of course, you need to show that the image is actually in $Y$,
namely, that the constraints that you put on the coordinates of $Y$
(typically something like the triangle inequality for the edges of a
triangle) are indeed satisfied for each geometric object of $X$.

This map $f$ will always be obviously continuous, and it is not too hard 
to show that it is proper and injective.  (Recall that a mapping is said 
to be proper if the pullback of a compact set is compact.) Then, the 
following lemma can be used to show that the image of $f$ is a union of 
connected components of $Y$.

\begin{lem}\label{PMT}
Let $X$ and $Y$ metric spaces, and let $f:X \rightarrow Y$ be a proper
local homeomorphism.  Then the image of $f$ is a union of connected 
components of $Y$.
\end{lem}

\noindent {\bf Proof of Lemma \ref{PMT}:}
It is sufficient to show that $f(X)$ is both open and closed in $Y$.  
Because $f$ is a local homeomorphism, it is an open mapping, so $f(X)$ is 
open in Y; and since $f: X \rightarrow Y$ is proper, it immediately 
follows that the limit of any sequence in the image of $f$ which converges 
in $Y$ must lie in the image of $f$, so $f(X)$ is closed in $Y$. 
\Endproof

In fact, a stronger result is true: any local homeomorphism between metric 
spaces which is also proper will be a finite-sheeted covering map \cite[p. 
23]{DOUADY} and \cite[p. 127]{LIMA}.  This gives an alternative route to 
proving Lemma \ref{PMT}.

Therefore, this lemma reduces the problem to showing that $X$ is nonempty 
and that $Y$ is connected, which are usually the hardest parts!

The result of the ``method of continuity'' is that you have established a
bijection between your geometric objects, set $X$, and the measurements $Y$.

\vspace{.1in}

Let $C$ be a cell complex on $\mathbb{S}^2$ that describes the 
combinatorics of a convex polyhedron.  We say that a hyperbolic polyhedron 
$P \subset \mathbb{H}^3$ {\em realizes} $C$ if there is a cellular 
homeomorphism from $C$ to $\partial P$ (i.e., a homeomorphism mapping 
faces of $C$ to faces of $P$, edges of $C$ to edges of $P$, and vertices 
of $C$ to vertices of $P$.) We will call each isotopy class of cellular 
homeomorphisms $\phi : C \rightarrow \partial P$ a {\em marking} on $P$.

Let $\Delta$ be the cell complex on $\mathbb{S}^2$ describing the 
combinatorics of the tetrahedron. Throughout this paper we will call 
hyperbolic polyhedra realizing $\Delta$ hyperbolic tetrahedra.

We will define ${\cal P}_\Delta$ to be the set of pairs $(P,\phi)$ so that 
$P$ is a hyperbolic tetrahedron and $\phi$ is a marking on $P$ with the 
equivalence relation that $(P,\phi) \sim (P',\phi ')$ if there exists an 
automorphism $\rho :  \mathbb{H}^3 \rightarrow \mathbb{H}^3$ such that 
$\rho(P) = P'$ and both $\phi '$ and $\rho \circ \phi$ represent the same 
marking on $P'$.

\begin{prop} \label{MANIFOLD}
The space ${\cal P}_\Delta$ is a manifold of dimension $6$.
\end{prop}

\noindent
{\bf Proof:}
Let ${\cal H}$ be the space of closed half-spaces of $\mathbb{H}^3$; 
clearly ${\cal H}$ is a 3-dimensional manifold.  Let ${\cal O}_\Delta$ be 
the set of marked hyperbolic polyhedra realizing $\Delta$. By forgetting 
this marking, an element of ${\cal O}_\Delta$ is a $4$-tuple of 
half-spaces that intersect in a polyhedron realizing $\Delta$.  This 
induces a mapping from ${\cal O}_\Delta$ to ${\cal H}^4$ whose image is an 
open set. We give ${\cal O}_\Delta$ the topology that makes this mapping 
from ${\cal O}_\Delta$ into ${\cal H}^4$ a local homeomorphism. Since 
${\cal H}^4$ is a $12$-dimensional manifold, ${\cal O}_\Delta$ must be a 
$12$-dimensional manifold as well.

If $\rho(P,\phi) = (P,\phi)$, we have that $\rho \circ \phi$ is isotopic 
to $\phi$ through cellular homeomorphisms.  Hence, the automorphism $\rho$ 
must fix all vertices of $P$, and consequently restricts to the identity 
on all edges and faces.  However, an automorphism of $\mathbb{H}^3$ which 
fixes four non-coplanar points must be the identity. Therefore ${\rm 
Aut}(\mathbb{H}^3)$ acts freely on ${\cal O}_\Delta$. This quotient is 
${\cal P}_\Delta$, hence ${\cal P}_\Delta$ is a manifold with dimension 
equal to ${\rm dim}({\cal O}_\Delta) - {\rm dim}({\rm Aut}(\mathbb{H}^3)) 
= 3 \cdot 4-6 = 6.$ 
\Endproof

In fact, we will restrict to the subset ${\cal P}_\Delta^0$ of tetrahedra 
with dihedral angles in $(0, \pi/2]$.  Notice that ${\cal P}_\Delta^0$ is 
not, {\it a priori}, a manifold or even a manifold with boundary.  All 
that we will need for the proof of Theorem \ref{TETRA1} is that ${\cal 
P}_\Delta$ is a manifold and that the subspace ${\cal P}_\Delta^0$ is a 
metric space.

Consider the map $\alpha: {\cal P}_\Delta \rightarrow \mathbb{R}^6$ which 
is obtained by measuring the dihedral angles (ordered by the marking) of 
an element of $P_\Delta$.  Using the topology on ${\cal P}_\Delta$ that is 
described in the proof of Proposition \ref{MANIFOLD}, it is clear that 
$\alpha$ is continuous. Therefore, we will use the method of continuity to 
show that $\alpha$ restricted to ${\cal P}_\Delta^0$ is a homeomorphism 
onto $A_\Delta$, in order to prove Theorem \ref{TETRA1}.

At this point it is necessary to clarify the statement of uniqueness in 
Theorem \ref{TETRA1}. We will show that the map $\alpha$ is injective, 
which shows that for each set of proposed dihedral angles 
$\alpha_1,\cdots,\alpha_6$ there is a unique {\it marked tetrahedron} with 
the dihedral angles $\alpha_1,\cdots,\alpha_6$, as ordered by this 
marking. This is what we mean by uniqueness in Theorem \ref{TETRA1} and in 
the later Theorem \ref{TETRA2}.

\vspace{.1in}

\noindent {\bf Proof of Theorem \ref{TETRA1}:}

The first step is to make sure that the dihedral angles of a compact 
tetrahedron satisfy conditions (1-3). For condition (1), notice that if 
two adjacent faces intersect along a line segment with dihedral angle $0$, 
they would coincide. In addition, the dihedral angle between adjacent 
faces is $\leq \pi/2$ by hypothesis. For condition (2), let $x$ be a 
vertex of $P$.  The compactness of $P$ implies that $x \in \mathbb{H}^3$, 
and by Lemma \ref{INTERSECT}, the sum of the dihedral angles between the 
three planes intersecting at $x$ must be $> \pi$. Furthermore, each face 
of a hyperbolic tetrahedron is a hyperbolic triangle of non-zero area so 
the Gauss-Bonnet formula gives condition (3). Therefore conditions (1-3) 
are necessary.

\vspace{.05in}

There is an elementary proof that $\alpha: {\cal P}_\Delta \rightarrow 
\mathbb{R}^E$ is injective: Since the face angles are uniquely determined 
by the dihedral angles and each face is a hyperbolic triangle, one can 
calculate the length of each edge using the hyperbolic law of cosines.

\vspace{.05in}

Before proving that $\alpha: {\cal P}_\Delta^0 \rightarrow A_\Delta$ is 
proper, we will need the following lemma:

\begin{lem} \label{NORMALIZE}
Given three points $v_1,v_2,v_3$ that form a non-obtuse, non-degenerate 
triangle in the Poincar\'e model of $\mathbb{H}^3$, there is a unique 
orientation preserving isometry taking $v_1$ to a positive point on the 
$x$-axis, $v_2$ to a positive point on the $y$-axis, and $v_3$ to a 
positive point on the $z$-axis. 
\end{lem}

\noindent {\bf Proof of Lemma \ref{NORMALIZE}:}
The points $v_1,v_2,$ and $v_3$ form a triangle $T$ in a plane $P_T$. It 
is sufficient to show that there is a plane $Q_T$ in the Poincar\'e ball 
model that intersects the positive octant in a triangle isomorphic to $T$.  
The isomorphism taking $v_1,v_2,$ and $v_3$ to the $x,y$, and $z$-axes 
will then be the one that takes the plane $P_T$ to the plane $Q_T$ and the 
triangle $T$ to the intersection of $Q_T$ with the positive octant.

Let $s_1,s_2,$ and $s_3$ be the side lengths of $T$.  The plane $Q_T$ must 
intersect the $x,y,$ and $z$-axes at distances $a_1,a_2,$ and $a_3$ 
satisfying the hyperbolic Pythagorean theorem:

$$\cosh(s_1) = \cosh(a_2) \cosh(a_3),$$
$$\cosh(s_2) = \cosh(a_3) \cosh(a_1),$$
$$\cosh(s_3) = \cosh(a_1) \cosh(a_2).$$

These equations can be solved for 
$(\cosh^2(a_1),\cosh^2(a_2),\cosh^2(a_3))$, obtaining $$\left( 
\frac{\cosh(s_2) \cosh(s_3)}{\cosh(s_1)}, \frac{\cosh(s_3) 
\cosh(s_1)}{\cosh(s_2)}, \frac{\cosh(s_1) \cosh(s_2)}{\cosh(s_3)} 
\right).$$ The only concern in solving for $a_i$ is that each of these 
terms is $\geq 1$. However, this follows from the triangle $T$ being 
non-obtuse. \Endproof

\begin{lem} \label{PROPER}
The mapping $\alpha: {\cal P}_\Delta^0 \rightarrow A_\Delta$ is proper.
\end{lem}

\noindent {\bf Proof:}

To see that $\alpha: {\cal P}_\Delta^0 \rightarrow A_\Delta$ is a proper 
mapping, suppose that there is a sequence of polyhedra $P_i$ realizing 
$\Delta$, with $\alpha(P_i) = {\bf a}_i \in A_\Delta$.  We must show that 
if ${\bf a}_i$ converges to ${\bf a} \in A_\Delta$, then a subsequence of 
the $P_i$ converges to some $P_\infty$ in ${\cal P}_\Delta^0$.

Throughout this part of the proof, we consider each $P_i$ to be in the 
Poincar\'e ball. Denote the vertices of $P_i$ by $v_1^i,v_2^i,v_3^i$, and 
$v_4^i$.  According to Lemma \ref{NORMALIZE}, we can normalize each $P_i$ 
so that $v_1^i$ is on the $x$-axis $v_2^i$ is on the $y$-axis, and $v_3^i$ 
is on the $z$-axis.

Because $\overline{\mathbb{H}^3}$ is a compact space (in the Euclidean 
metric), we can take a subsequence of the $P_i$ so that the vertices 
$v_1^i,\cdots,v_4^i$ converge to some points $v_1,\cdots,v_4$ in 
$\overline{\mathbb{H}^3}$.  We must use that ${\bf a}$ satisfies 
conditions (1-3) to show that $v_1,\cdots,v_4$ are actually at distinct 
finite points in $\mathbb{H}^3$ whose span is a tetrahedron.

\vspace{.05in}
\noindent
{\bf The vertices $v_1^i, v_2^i, v_3^i$, and $v_4^i$ converge to
distinct points in $\overline{\mathbb{H}^3}$}
\vspace{.05in}

Notice that at most two of the vertices could converge to the same point 
in $\partial \mathbb{H}^3$, since $v_1^i$ is on the $x$-axis, $v_2^i$ is 
on the $y$-axis, and $v_3^i$ is on the $z$-axis. We suppose, without loss 
of generality, that $v_4^i$ converges to the same point in $\partial 
\mathbb{H}^3$ as $v_3^i$, that is, both $v_4^i$ and $v_3^i$ converge to 
the north pole of the Poincar\'e ball.  Then, however, the dihedral angle, 
$\psi$, between the face spanned by $(v_1^i, v_2^i, v_3^i)$ and the face 
spanned by $(v_1^i, v_2^i, v_4^i)$ must limit to $0$, contrary to 
condition (1).  This configuration is depicted in the diagram below.

\vspace{.1in}
\begin{center}
\begin{picture}(0,0)%
\includegraphics{tetra_nc.pstex}%
\end{picture}%
\setlength{\unitlength}{4144sp}%
\begingroup\makeatletter\ifx\SetFigFont\undefined%
\gdef\SetFigFont#1#2#3#4#5{%
  \reset@font\fontsize{#1}{#2pt}%
  \fontfamily{#3}\fontseries{#4}\fontshape{#5}%
  \selectfont}%
\fi\endgroup%
\begin{picture}(2853,2847)(580,-2581)
\put(2747,-1078){\makebox(0,0)[lb]{\smash{\SetFigFont{8}{9.6}{\familydefault}{\mddefault}{\updefault}{\color[rgb]{0,0,0}\small{$v_2^i$}}%
}}}
\put(1257,-1358){\makebox(0,0)[lb]{\smash{\SetFigFont{8}{9.6}{\familydefault}{\mddefault}{\updefault}{\color[rgb]{0,0,0}\small{$v_1^i$}}%
}}}
\put(2099, 31){\makebox(0,0)[lb]{\smash{\SetFigFont{8}{9.6}{\familydefault}{\mddefault}{\updefault}{\color[rgb]{0,0,0}\small{$v_4^i$}}%
}}}
\put(1785, 31){\makebox(0,0)[lb]{\smash{\SetFigFont{8}{9.6}{\familydefault}{\mddefault}{\updefault}{\color[rgb]{0,0,0}\small{$v_3^i$}}%
}}}
\put(2113,-1378){\makebox(0,0)[lb]{\smash{\SetFigFont{8}{9.6}{\familydefault}{\mddefault}{\updefault}{\color[rgb]{0,0,0}\small{$\psi$}}%
}}}
\end{picture}

\end{center}
\vspace{.1in}

Hence, we conclude that any of the vertices $v_j^i$ that converge to 
points in $\partial \mathbb{H}^3$, must converge to distinct points.

Any face of $P_i$ that degenerates to a point or a line segment has 
(hyperbolic) area that limits to zero, since the vertices of $P_i$ that 
converge to points in $\partial \mathbb{H}^3$ converge to distinct points. 
Hence, by the Gauss Bonnet formula, the sum of the face angles for such a 
degenerating face would limit to $\pi$, contrary to condition (3). This is 
enough to show that $v^i_1,\cdots,v^i_4$ converge to distinct points 
$v_1,\cdots,v_4$ in $\overline{\mathbb{H}^3}$.

\vspace{.05in}
\noindent
{\bf The limit points $v_1, v_2, v_3$, and $v_4$ are finite
points whose span is a tetrahedron.}
\vspace{.05in}

The sum of the dihedral angles at the edges leading to each $v_j^i$ 
converges to a value $> \pi$.  Therefore, according to Lemma 
\ref{INTERSECT}, we conclude that the limit points of vertices 
$v_1,\cdots,v_4$ are actually at finite points.

Since each face is non-degenerate, and the dihedral angles are non-obtuse, 
the $P_i$ cannot degenerate to a single triangle.  So, their span realizes 
a tetrahedron, with dihedral angles ${\bf a}$.

This is enough to conclude that 
$\alpha: {\cal P}_\Delta^0 \rightarrow A_\Delta$ is proper.
\Endproof

\vspace{.05in}

Invariance of Domain gives that $\alpha: {\cal P}_\Delta \rightarrow 
\mathbb{R}^6$ is a local homeomorphism because it is a continuous and 
injective mapping between manifolds of the same dimension. Therefore, the 
restriction $\alpha: {\cal P}^0_\Delta \rightarrow A_\Delta$ is also a 
local homeomorphism.  Because $\alpha: {\cal P}^0_\Delta \rightarrow 
A_\Delta$ is also a proper mapping, by Lemma \ref{PMT}, $\alpha({\cal 
P}_\Delta^0)$ is a union of connected components of $A_\Delta$. We will 
show that ${\cal P}^0_\Delta$ is nonempty and that $A_\Delta$ is 
connected, thus proving that $\alpha: {\cal P}^0_\Delta \rightarrow 
A_\Delta$ is surjective.

\vspace{.1in}
                                                                                
The easiest way to see that ${\cal P}^0_\Delta \neq \emptyset$ is by 
explicit construction. Let ${\bf v}_1 = (0,1,0,0), {\bf v}_2 = 
(0,0,1,0),{\bf v}_3 = (0,0,0,1),$ and ${\bf v}_4 = \frac{1}{\sqrt 
2}(-1,-1,-1,-1)$.  Then the intersection of the half-spaces $H_{{\bf v}_1} 
\cap H_{{\bf v}_2} \cap H_{{\bf v}_3} \cap H_{{\bf v}_4}$ is a hyperbolic 
tetrahedron with dihedral angles $\alpha_{1,2} = \pi/2, \alpha_{1,3} = 
\pi/2, \alpha_{2,3} = \pi/2, \alpha_{1,4} = \alpha_{2,4} = \alpha_{3,4} = 
\arccos(1/\sqrt 2) = \pi/4$.  Hence, we conclude that ${\cal P}^0_\Delta 
\neq \emptyset.$

To see that $A_\Delta$ is connected is significantly harder than for $A_C$ 
with $C$ not the tetrahedron because the inequalities specifying 
$A_\Delta$ are not linear. We will have to do detailed analysis of the 
equation that expresses a face's angles in terms of the dihedral angles.
                                                                                
\begin{lem}
$A_\Delta$ is path connected.
\end{lem}
                                                                                
\noindent {\bf Proof:}
Recall from Lemma \ref{TRIVALENT} that the face angle $\beta_i$ at a vertex $(e_i,e_j,e_k)$ in the face
containing $e_j$ and $e_k$ is 
                                                                                
$$\cos(\beta_i) = \frac{\cos(\alpha_i) +\cos(\alpha_j)\cos(\alpha_k)}
{\sin(\alpha_j)\sin(\alpha_k)}$$

Let $A_i \subset \partial A_\Delta$ be the subset obtained by restricting 
the dihedral angle sum at each of the vertices, except $v_i$, to equal 
$\pi$.  Using the formula for the $\beta_j$, one can check that at each 
vertex with dihedral angle sum exactly $\pi$, each of the face angles is 
$0$.  One can also check that each of the face angles at $v_i$ is 
non-obtuse, since each of the dihedral angles is non-obtuse.  Therefore, 
for any point in $A_i$, for each $i=1,\cdots,4$, each of the face angle 
sums is $\leq \pi/2$.  Therefore, since the formula for face angles in 
terms of dihedral angles is continuous, there exists a neighborhood $NA_i$ 
of each $A_i$ in $A_\Delta$.  If necessary, we can restrict $NA_i$ to a 
smaller set which is connected, since $A_i$ is convex.

For $i=1,\cdots,4$, each $A_i$ contains $(\pi/3,\cdots,\pi/3)$, which are 
the dihedral angles of the regular ideal tetrahedron, hence $NA_1 \cap 
NA_2 \cap NA_3 \cap NA_4 \neq \emptyset$.  Therefore $NA_1 \cup NA_2 \cup 
NA_3 \cup NA_4$ is path connected.  Denote this set by $\cal{N}.$

\vspace{.05in}
                                                                                
Given any ${\bf a} \in A_\Delta$, we will create a path from ${\bf a}$ to 
a point in $\cal{N}$.  This will be sufficient to prove that $A_\Delta$ is 
connected.  First, notice that for any ${\bf a} \in A_\Delta$, decreasing 
any of the components of ${\bf a}$ does not increase any of the $\beta_i$.  
One can check that if:
                                                                                
$$F(x,y,z) = \frac{\cos(x)+\cos(y)\cos(z)}{\sin(y)\sin(z)}$$
                                                                                
\noindent
Then we have:
                                                                                
$$\frac{\partial F}{\partial x} = -\frac{\sin(x)}{\sin(y)\sin(z)},$$
$$\frac{\partial F}{\partial y} = \frac{-\sin(y)\sin(z)\sin(y)\cos(z) -
\cos(y)\cos(z)\cos(y)\sin(z)}{\sin^2(y)\sin^2(z)},$$
$$\frac{\partial F}{\partial z} = \frac{-\sin(y)\sin(z)\cos(y)\sin(z) -
\cos(y)\cos(z)\sin(y)\cos(z)}{\sin^2(y)\sin^2(z)}.$$
                                                                                
These have the nice property that for all $(x,y,z) \in (0,\pi/2]^3$ we
have $\frac{\partial F}{\partial x} < 0, \frac{\partial F}{\partial y} <
0$, and $\frac{\partial F}{\partial z} < 0$.  Because $\arccos$ is a
decreasing function, this gives that $\beta(\gamma_i,\gamma_j,\gamma_k)
\leq \beta({\bf a}_i,{\bf a}_j,{\bf a}_k)$ when $\gamma_i \leq {\bf a}_i,
\gamma_j \leq {\bf a}_j,$ and $\gamma_k \leq {\bf a}_k$.
                                                                                
\noindent
Therefore, given ${\bf a} \in A_\Delta$, decreasing the angles of ${\bf a}$
cannot result in a violation of condition (3).
                                                                               
Consider $t \cdot {\bf a}$ decreasing $t$ from $1$ to $0$.  For some first 
value of $t$, the sum of dihedral angles at one of the vertices, say 
$v_1$, will be $\pi$.  Next, decrease only the dihedral angles of edges 
not entering $v_1$ in the same uniform way until the sum of the dihedral 
angles at another of the vertices, say $v_2$ equals $\pi$.  Finally, 
decrease the dihedral angle on the edge that does not enter $v_1$ or $v_2$ 
until one the two remaining vertices has dihedral angle sum $\pi$, call 
this vertex $v_3$.
                                                                                
Since we have decreased the dihedral angles during the duration of this
path, condition (3) was satisfied throughout.  Condition (1) was satisfied
throughout because we decreased the dihedral angles, so none exceeded
$\pi/2$ and since we decreased them by scaling, so that none reached $0$.

Therefore, we have constructed a path from ${\bf a}$ to $A_1$. This path 
must have entered $\cal{N}$ because it connected the point ${\bf a} \in A$ 
to $A_1$.
\Endproof
                                                                                
Therefore, since $\alpha_\Delta : {\cal P}^0_\Delta \rightarrow A_\Delta$ is
an injective covering map with ${\cal P}^0_\Delta \neq \emptyset$
and $A_\Delta$ path connected, we conclude that $\alpha_\Delta$ is a
homeomorphism.  This proves Theorem \ref{TETRA1}. \Endproof

Using equation \ref{TLC}, we can re-express
Theorem \ref{TETRA1} entirely in terms of the dihedral angles.

\begin{cor} \label{TETRA2}
                                                                                
There is a compact hyperbolic tetrahedron
with non-obtuse dihedral angles
$\alpha_1,\cdots,\alpha_6$ if and only if:
                                                                                
\begin{enumerate}
                                                                                
\item  For each edge $e_i$, $0 < \alpha_i \leq \pi/2$.
                                                                                
\item  Whenever 3 (distinct) edges $e_i, e_j, e_k$ meet at a vertex,
$\alpha_i + \alpha_j+\alpha_k > \pi$.
                                                                                
\item  For each face $F$ bounded by edges $e_i,e_j,e_k$ with edges
$e_{i,j},e_{j,k},e_{k,i}$ emanating from the vertices, we have:
                                                                                
{\small
$$
\arccos\left(\frac{\cos(\alpha_{i,j}) + \cos(\alpha_i)\cos(\alpha_j)}{\sin(\alpha_i)\sin(\alpha_j)}\right) + \hspace{.4in}$$
                                                                                
$$\arccos\left(\frac{\cos(\alpha_{j,k}) + \cos(\alpha_j)\cos(\alpha_k)}{\sin(\alpha_j)\sin(\alpha_k)}\right) + \hspace{.3in} $$
                                                                                
$$\arccos\left(\frac{\cos(\alpha_{k,i}) + \cos(\alpha_k)\cos(\alpha_i)}{\sin(\alpha_k)\sin(\alpha_i)}\right) < \pi .$$
}
                                                                                
\end{enumerate}

Furthermore, this hyperbolic polyhedron is unique.
                                                                                
\end{cor}
                                                                                
\noindent The proof is evidently a direct consequence of Theorem \ref{TETRA1}
and the formula for the face angles.
                                                                                
\vspace{.1in}

The reader should notice that Theorem \ref{TETRA1} can also be proved directly
from Theorem \ref{GRAM}, the characterization of hyperbolic tetrahedra in terms
of the Gram matrix $M$.  In Lemma \ref{INTERSECT} we checked that if we
restrict to non-obtuse dihedral angles, Condition (2) from Theorem \ref{TETRA1}
is equivalent to the condition that every principal minor of $M$ is positive
definite.

Similar trigonometric tricks can be used to show that Conditions (1) and (3)
from Theorem \ref{TETRA1} are equivalent to $\det(M) < 0$.  If a face $F$
contains face angles $\beta_i, \beta_j,\beta_k$ and the edges surrounding $F$
have dihedral angles $\alpha_l,\alpha_m,\alpha_n$, then: 
{\footnotesize
\begin{eqnarray*} && \det(M) = \\ &&
-4(1-\cos^2\alpha_l)(1-\cos^2\alpha_m)(1-\cos^2\alpha_n)\cdot \\ &&
\cos\left(\frac{\beta_i+\beta_j+\beta_k}{2}\right)\cos\left(\frac{\beta_i-\beta_j+\beta_k}{2}\right)\cos\left(\frac{\beta_i+\beta_j-\beta_k}{2}\right)\cos\left(\frac{-\beta_i+\beta_j+\beta_k}{2}\right).
\end{eqnarray*} }

Condition (1) from Theorem \ref{TETRA1} requires that the dihedral angles are
positive and non-obtuse, so the second line of this equation is negative.  The
third line is positive if and only if $\beta_i+\beta_j+\beta_k < \pi$, since
the face angles are non-obtuse (because the dihedral angles are non-obtuse.)
Hence, Theorem \ref{TETRA1} does follow from Theorem \ref{GRAM}.  However, the
author feels that the proof using the method of continuity is more intuitive.

\vspace{.1in}
In terms of the dihedral angles, condition (3) is reasonably nasty.  In 
fact, it results in $A_\Delta$ being non-convex!  Consider the hyperbolic 
tetrahedron with dihedral angles $x$ and $y$ assigned to two edges that 
meet at a vertex and dihedral angles $\alpha$ assigned to the remaining 
edges. The following figure was computed in Maple \cite{MAP} and shows the 
cross section of $A_\Delta$ when $\alpha = 1.3$.
                                                                                
\vspace{.1in}
\begin{center}
\includegraphics[scale=.45]{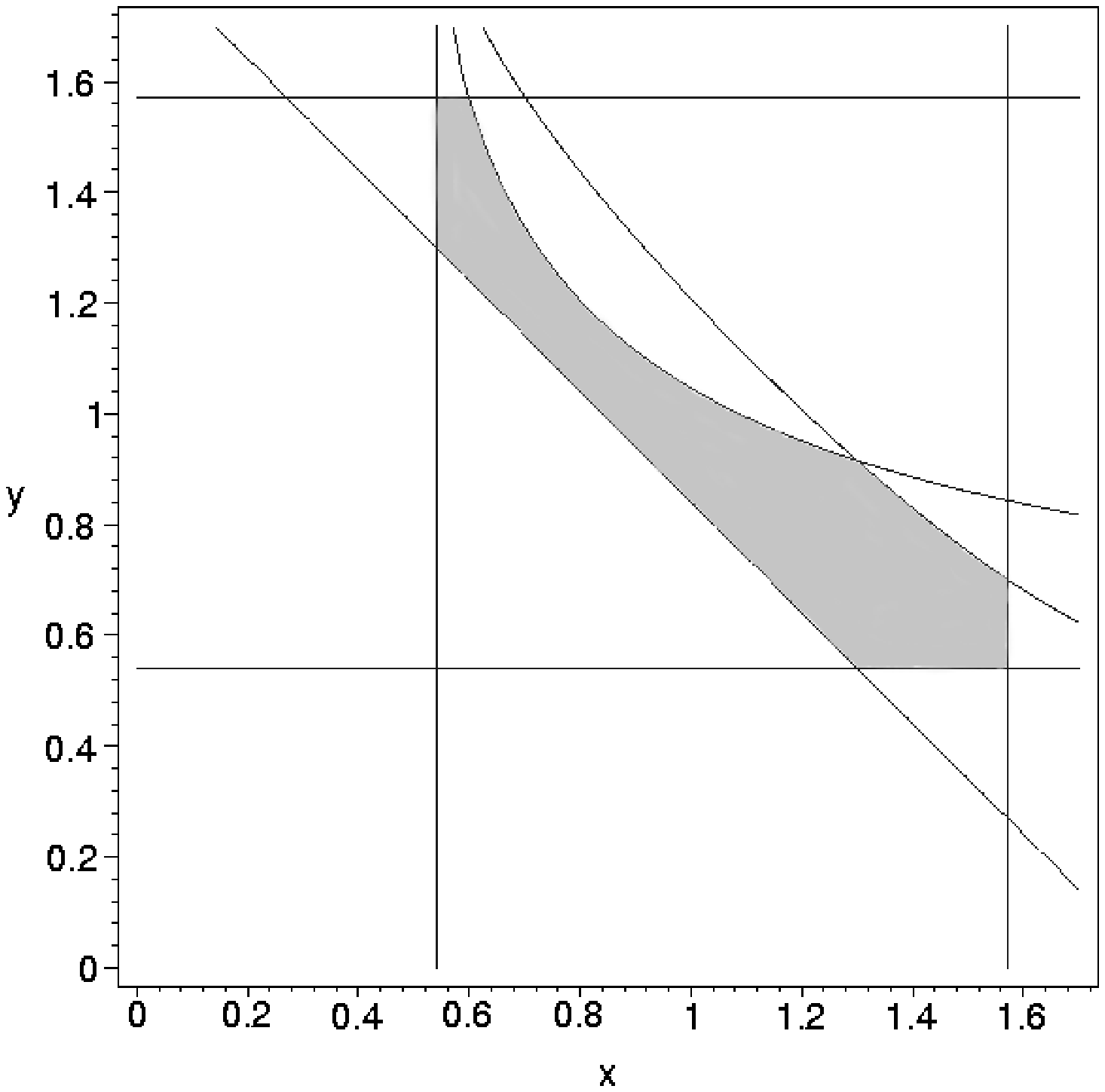}
\end{center}
                                                                                
\vspace{.1in}

This classification of hyperbolic tetrahedra in terms of their dihedral 
angles gives us some understanding of how a generalization of Andreev's 
Theorem \cite{AND,H,ROE2,ROE} to include obtuse dihedral angles would be 
significantly more complicated than Andreev's Theorem.

One obvious difficulty in considering arbitrary dihedral angles is that 
one cannot restrict to studying hyperbolic polyhedra realizing trivalent 
abstract polyhedra, a restriction that was essential in the proof of 
Andreev's Theorem.

However, a further difficulty that arises even for trivalent hyperbolic 
polyhedra is that for each $n$-sided face ${\cal F}$, there is the 
necessary condition that the sum of the face angles of ${\cal F}$ must be 
$< (n-2)\pi$, resulting from the Gauss-Bonnet Theorem. As in conditions 
(3) in Theorems \ref{TETRA1} and \ref{TETRA2} from this paper, this 
results in highly non-linear necessary conditions on the dihedral angles.

The restriction to non-obtuse dihedral angles results in non-obtuse face 
angles via Equation (\ref{TLC}).  Therefore, when studying hyperbolic 
polyhedra with non-obtuse dihedral angles, this condition on face angles 
is irrelevant for faces with more than $5$ edges. Part of the proof of 
Andreev's Theorem is to show that, as long as the polyhedron has more than 
4 faces, this condition on face angles for 3 and 4-sided faces is 
automatically a consequence of two other {\em linear} necessary conditions 
on the dihedral angles. In the statement of Andreev's Theorem as written 
in \cite{ROE2,ROE}, these are conditions (3-5).

However, once the dihedral angles are non-obtuse, these conditions on face 
angles for 3 and 4-sided faces are no longer a consequence of conditions 
(3-5) of Andreev's Theorem.  Furthermore, this condition on face angles 
becomes relevant for faces with 5 and more edges because the face angles 
are no longer restricted to be non-obtuse.

Of course, one can also expect that other conditions may be necessary to 
prevent more exotic types of degeneracies.

\vspace{.2in}
\noindent
{\bf \Large Acknowledgments}

\noindent 
The author thanks the anonymous referee for his or her detailed
comments, and especially for informing the author of the characterization of
hyperbolic tetrahedra in terms of the Gram matrix.  The author also
thanks Joshua Bowman, William D. Dunbar, John H. Hubbard, and Rodrigo Perez
for their mathematical suggestions.

\bibliographystyle{plain}
\bibliography{andreev.bib}
\end{document}